\documentclass[11pt]{amsart}
\usepackage{latexsym,graphicx}

\numberwithin{equation}{section}
\theoremstyle{plain}

\theoremstyle{remark}

\theoremstyle{definition}

\newcommand{\D}{{\mathcal D}}
\newcommand{\E}{\mathcal E}

\newcommand{\G}{{\mathcal G}}

\newcommand{\K}{{\mathcal K}}
\renewcommand{\L}{{\mathcal L}}
\newcommand{\M}{{\mathcal M}}
\newcommand{\N}{\mathbb N}

\newcommand{\V}{{\mathcal V}}

\newcommand{\dist}{\operatorname{dist}}

\newcommand{\fp}{\operatorname{FP}}

\newcommand{\Int}{\operatorname{Int}}

\newcommand{\on}{\:\mbox{\rule{0.1ex}{1.2ex}\rule{1.9ex}{0.1ex}}\:}

\renewcommand{\span}{\operatorname{span}}

\newcommand{\supp}{\operatorname{Supp}}

\def\Ga{\Gamma}
\def\half{{1 \over 2}}

\newcommand{\oa}{\overrightarrow}

\newcommand{\ol}{\overline}

\def\XXint#1#2#3{{\setbox0=\hbox{$#1{#2#3}{\int}$}
      \vcenter{\hbox{$#2#3$}}\kern-.5\wd0}}

\begin{document}

\def\cal{\mathcal}

\font\tpt=cmr10 at 12 pt
\font\fpt=cmr10 at 14 pt

\font \fr = eufm10




\overfullrule=0in

\def\boxit#1{\hbox{\vrule
 \vtop{%
  \vbox{\hrule\kern 2pt %
     \hbox{\kern 2pt #1\kern 2pt}}%
   \kern 2pt \hrule }%
  \vrule}}

  \def\harr#1#2{\ \smash{\mathop{\hbox to .3in{\rightarrowfill}}\limits^{\scriptstyle#1}_{\scriptstyle#2}}\ }

\def\ALL{1}
\def\BTA{2}
\def\BL{3}
\def\BRE{4}
\def\CNS{5}
\def\CIL{6}
\def\CRA{7}
\def\DDD{8}
\def\DDR{9}
\def\GEO{10}
\def\HYP{11}
\def\BEL{12}
\def\AC{13}
\def\SURVEY{14}
\def\NOTES{15}
\def\AET{16}
\def\LAG{17}
\def\KRY{18}
\def\PLI{19}
\def\RT{20}
\def\SLO{21}
\def\TRU{22}
\def\TWC{23}
\def\WAL{24}

 \def\GG{{{\bf G} \!\!\!\! {\rm l}}\ }

\def\GL{{\rm GL}}

\def\bll{I \!\! L}

\def\IFF{\qquad\iff\qquad}
\def\bra#1#2{\langle #1, #2\rangle}
\def\bbf{{\bf F}}
\def\bbj{{\bf J}}
\def\Jtn{{\bbj}^2_n}  \def\JtN{{\bbj}^2_N}  \def\JoN{{\bbj}^1_N}
\def\jt{j^2}
\def\jtx{\jt_x}
\def\Jt{J^2}
\def\Jtx{\Jt_x}
\def\bpp{{\bf P}^+}
\def\bpt{{\wt{\bf P}}}
\def\fsh{$F$-subharmonic }
\def\mo{monotonicity }
\def\jet{(r,p,A)}
\def\ss{\subset}
\def\sse{\subseteq}
\def\half{\hbox{${1\over 2}$}}
\def\smfrac#1#2{\hbox{${#1\over #2}$}}
\def\oa#1{\overrightarrow #1}
\def\dim{{\rm dim}}
\def\dist{{\rm dist}}
\def\codim{{\rm codim}}
\def\deg{{\rm deg}}
\def\rank{{\rm rank}}
\def\log{{\rm log}}
\def\Hess{{\rm Hess}}
\def\Hessyp{{\rm Hess}_{\rm SYP}}
\def\trace{{\rm trace}}
\def\tr{{\rm tr}}
\def\max{{\rm max}}
\def\min{{\rm min}}
\def\span{{\rm span\,}}
\def\Hom{{\rm Hom\,}}
\def\det{{\rm det}}
\def\End{{\rm End}}
\def\Sym{{\rm Sym}^2}
\def\diag{{\rm diag}}
\def\pt{{\rm pt}}
\def\Spec{{\rm Spec}}
\def\pr{{\rm pr}}
\def\Id{{\rm Id}}
\def\Grass{{\rm Grass}}
\def\Herm#1{{\rm Herm}_{#1}(V)}
\def\arr{\longrightarrow}
\def\supp{{\rm supp}}
\def\Link{{\rm Link}}
\def\Wind{{\rm Wind}}
\def\Div{{\rm Div}}
\def\vol{{\rm vol}}
\def\foral{\qquad {\rm for\ all\ \ }}
\def\fpsh{{\cal PSH}(X,\f)}
\def\Core{{\rm Core}}
\def\dis{f_M}
\def\Re{{\rm Re}}
\def\rn{\bbr^n}
\def\pp{\cp^+}
\def\plp{\cp_+}
\def\Int{{\rm Int}}
\def\cix{C^{\infty}(X)}
\def\Gr#1{G(#1,\rn)}
\def\Symn{{\Sym(\rn)}}
\def\SymN{{\Sym(\bbr^N)}}
\def\Gpn{G(p,\rn)}
\def\fd{{\rm free-dim}}
\def\SA{{\rm SA}}
 \def\cd{{\cal C}}
 \def\cdt{{\widetilde \cd}}
 \def\cm{{\cal M}}
 \def\cmt{{\widetilde \cm}}

\def\Theorem#1{\medskip\noindent {\bf THEOREM \bf #1.}}
\def\Prop#1{\medskip\noindent {\bf Proposition #1.}}
\def\Cor#1{\medskip\noindent {\bf Corollary #1.}}
\def\Lemma#1{\medskip\noindent {\bf Lemma #1.}}
\def\Remark#1{\medskip\noindent {\bf Remark #1.}}
\def\Note#1{\medskip\noindent {\bf Note #1.}}
\def\Def#1{\medskip\noindent {\bf Definition #1.}}
\def\Claim#1{\medskip\noindent {\bf Claim #1.}}
\def\Conj#1{\medskip\noindent {\bf Conjecture \bf    #1.}}
\def\Ex#1{\medskip\noindent {\bf Example \bf    #1.}}
\def\Qu#1{\medskip\noindent {\bf Question \bf    #1.}}
\def\Exercise#1{\medskip\noindent {\bf Exercise \bf    #1.}}

\def\HoQu#1{ {\AAA T\BBB HE\ \AAA H\BBB ODGE\ \AAA Q\BBB UESTION \bf    #1.}}

\def\pf{\medskip\noindent {\bf Proof.}\ }
\def\qed{\hfill  $\vrule width5pt height5pt depth0pt$}
\def\equdef{\buildrel {\rm def} \over  =}
\def\qedqed{\hfill  $\vrule width5pt height5pt depth0pt$ $\vrule width5pt height5pt depth0pt$}
\def\mathqed{  \vrule width5pt height5pt depth0pt}

\def\V{W}

\def\df{d^{\phi}}
\def\hk{\_{\rm l}\,}
\def\n{\nabla}
\def\w{\wedge}

\def\cu{{\cal U}}   \def\cc{{\cal C}}   \def\cb{{\cal B}}  \def\cz{{\cal Z}}
\def\cv{{\cal V}}   \def\cp{{\cal P}}   \def\ca{{\cal A}}
\def\cw{{\cal W}}   \def\co{{\cal O}}
\def\ce{{\cal E}}   \def\ck{{\cal K}}
\def\ch{{\cal H}}   \def\cm{{\cal M}}
\def\cs{{\cal S}}   \def\cn{{\cal N}}
\def\cd{{\cal D}}
\def\cl{{\cal L}}
\def\cp{{\cal P}}
\def\cf{{\cal F}}
\def\ccr{{\cal  R}}

\def\gerG{{\fr{\hbox{g}}}}
\def\gerB{{\fr{\hbox{B}}}}
\def\gerR{{\fr{\hbox{R}}}}
\def\p#1{{\bf P}^{#1}}
\def\vf{\varphi}

\def\wt{\widetilde}
\def\wh{\widehat}

\def\and{\qquad {\rm and} \qquad}
\def\arr{\longrightarrow}
\def\ol{\overline}
\def\bbr{{\mathbb R}}\def\bbh{{\mathbb H}}\def\bbo{{\mathbb O}}
\def\bbc{{\mathbb C}}
\def\bbq{{\mathbb Q}}
\def\bbz{{\mathbb Z}}
\def\bbp{{\mathbb P}}
\def\bbd{{\mathbb D}}

\def\a{\alpha}
\def\b{\beta}
\def\d{\delta}
\def\e{\epsilon}
\def\f{\phi}
\def\g{\gamma}
\def\k{\kappa}
\def\l{\lambda}
\def\o{\omega}

\def\s{\sigma}
\def\x{\xi}
\def\z{\zeta}

\def\D{\Delta}
\def\L{\Lambda}
\def\G{\Gamma}
\def\O{\Omega}

\def\bd{\partial}
\def\bdf{\partial_{\f}}
\def\lag{Lagrangian}
\def\psh{plurisubharmonic }
\def\ph{pluriharmonic }
\def\pph{partially pluriharmonic }
\def\omp{$\omega$-plurisubharmonic \ }
\def\ffl{$\f$-flat}
\def\PH#1{\widehat {#1}}
\def\lloc{L^1_{\rm loc}}
\def\dbar{\ol{\partial}}
\def\lp{\Lambda_+(\f)}
\def\lpp{\Lambda^+(\f)}
\def\bo{\partial \Omega}
\def\Ob{\overline{\O}}
\def\fc{$\phi$-convex }
\def\PSH{{ \rm PSH}}
\def\SH{{\rm SH}}
\def\totr{ $\phi$-free }
\def\BM{\lambda}
\def\Der{D}
\def\CH{{\cal H}}
\def\RH{\overline{\ch}^\f }
\def\pconv{$p$-convex}
\def\MA{MA}
\def\lagpsh{Lagrangian plurisubharmonic}
\def\hermsk{{\rm Herm}_{\rm skew}}
\def\PSHl{\PSH_{\rm Lag}}
 \def\ppsh{$\pp$-plurisubharmonic}
\def\fp{$\pp$-plurisubharmonic }
\def\fh{$\pp$-pluriharmonic }
\def\Symn{\Sym(\rn)}
 \def\ci{C^{\infty}}
\def\USC{{\rm USC}}
\def\LSC{{\rm LSC}}
\def\fa{{\rm\ \  for\ all\ }}
\def\ppc{$\pp$-convex}
\def\cpt{\wt{\cp}}
\def\ft{\wt F}
\def\ob{\overline{\O}}
\def\Be{B_\e}
\def\K{{\rm K}}

\def\M{{\bf M}}
\def\N#1{C_{#1}}
\def\ds{Dirichlet set }
\def\dir{Dirichlet }
\def\Fa{{\oa F}}
\def\TR{{\cal T}}
 \def\ISO{{\rm ISO_p}}
 \def\Span{{\rm Span}}

\def\ALL{1}
\def\AV{2}
\def\BTA{3}
\def\BL{4}
\def\BRE{5}
\def\CNS{6}
\def\CP{7}
\def\CPW{8}
\def\CIL{9}
\def\CRA{10}
\def\DTT{11}
\def\DON{12}
\def\CG{13}
\def\DDD{14}
\def\DDR{15}
\def\GEO{16}
\def\HYP{17}
\def\BEL{18}
\def\SURVEY{19}
\def\AC{20}
\def\NOTES{21}
\def\AET{22}
\def\TANG{23}
\def\TANGG{24}
\def\LAG{25}
\def\SLE{26}
\def\JTY{27}
\def\KRY{28}
\def\PLI{29}
\def\RT{30}
\def\SLO{31}
\def\SPR{32}
\def\TRUU{33}
\def\TRU{34}
\def\TWC{35}
\def\TWCC{36}
\def\TWCCC{37}
\def\WAL{38}

\vskip .4in

\def\E{E}
\def\fpsi{{F_f(\psi)}}
\def\bL{{\bf \Lambda}}
\def\bdf{{\bf f}}
\def\UU{U}
\def\bbm{{\bf M}}
\def\gg{{\mathfrak g}}
\def\iv{^{-1}}
\def\Sn{{\mathcal S}(n)}
\def\Ga{G\aa rding}
\def\GD{\Ga-Dirichlet}
\def\D{\mathbb D}
\def\gra{\d}
\def\on{^{1\over N}}

\font\headfont=cmr10 at 14 pt
\font\aufont=cmr10 at 11 pt

\def\AO{[AO]}
\def\ADO{[ADO]}
\def\Gar{[Ga]}
\def\GP{[GP$_1$]} 
\def\GGP{[GP$_2$]}
\def\GGGP{[GP$_3$]}
\def\GGGGP{[GP$_4$]}
\def\GPT{[GPT]} 
 \def\Gu{[Gu]}
\def\HL{[HL$_1$]}
\def\HHL{[HL$_2$]}
\def\HHLL{[HL$_3$]}
\def\HHHL{[HL$_4$]}
\def\HHHHL{[HL$_5$]}
\def\KLW{[KLW]}
\def\Ren{[Ren]}

\title[ DETERMINANT MAJORIZATION FOR  NONLINEAR OPERATORS]
{\headfont A  DEFINITIVE DETERMINANT MAJORIZATION RESULT FOR NONLINEAR OPERATORS}

\date{\today}
\author{ \aufont F. Reese Harvey and H. Blaine Lawson, Jr.}
\thanks
{Second author was partially supported by the Simons Foundation}

\maketitle


\centerline{\sl  Dedicated to John Polking}
\centerline{\sl on the occasion of his ninetieth birthday   }
\vskip .5in

\centerline{\bf Abstract}
Let $\gg$ be a \Ga-Dirichlet operator on the set $\Sn$ of symmetric $n\times n$ matrices.
We assume that $\gg$ is $I$-central, that is, $D_I \gg = k I$ for some $k>0$.  Then
$$
\gg(A)^{1\over N} \ \geq\ \gg(I)^{1\over N} (\det\, A)^{1\over n} \qquad \forall\, A>0.
$$

\noindent
From work of Guo, Phong, Tong,     Abja, Dinew, Olive and many others, this inequality has important applications.
\medskip

\vskip .5in

 \vfill \eject







\centerline
{\bf 1.\  Introduction.}  

This paper is concerned with G\aa rding-Dirichlet operators.  These are nonlinear elliptic differential operators 
on euclidean space, and in many important cases on riemannian manifolds, hermitian manifolds, etc.
(This is discussed in detail in \HL.)
These operators are defined by applying the following polynomials to the Hessians of functions.
Let $\Sn$ denote the space  of second derivatives, i.e., the vector space  of $n\times n$ real symmetric matrices.

\Def{1.1} Suppose $\gg$ is a real homogeneous polynomial of degree $N$ on  $\Sn$  which is $I$-hyperbolic,
which means that $\gg(I)>0$ and for each $A \in \Sn$, the polynomial $\vf_A(t)\equiv \gg(tI+A)$ has only real roots.
By the {\bf G\aa rding cone} $\G$ of $\gg$ we mean the connected component of $\{\gg>0\}$
containing $I$.  If $\G$ contains $\Int \cp = \{A:A>0\}$, then $\gg$ is called a {\bf \Ga-Dirichlet polynomial} 
or a {\bf \Ga-Dirichlet polynomial operator}. 

In \HHHL\ we proved a determinant majorization formula for such polynomials provided
they had certain invariance properties.   This had applications to work of B. Guo, D. Phong, F. Tong \,  \GPT, \GP,  \GGP, \GGGP\
 and S. Abja, S. Dinew, G. Olive  \AO, \ADO.  In fact, the important paper \GPT\ has lead to much subsequent work.
 The reader could check  recent papers \GGGGP\ and \KLW\ for history and references.

The point of this paper is to prove  a more  definitive result that the majorization  formula holds for a wide class 
of \Ga-Dirichlet polynomials where the invariance property is not needed.  One needs only to 
check the Central Ray Hypothesis, which is generally  very easy to do. (See Section 3.)

\Def{1.2} The {\bf Central Ray Hypothesis}\footnote{This concept was first introduced in \HHHL\
 with several equivalent formulations in Theorem E.8.}  is that:

{\sl The gradient of $\gg$  at $I$ satisfies $D_I \gg = k I$    for some $k>0$.}\footnote{In this paper $D_AF$ denotes the gradient of $F$ at the point $A$}  

\noindent
This is equivalent to the following:

{\sl The $\gg$-Laplacian $\Delta^\gg$, defined to be the sum of the $I$-eigenvalues of $\gg$,
is a positive constant times the standard Laplacian $\Delta^{\rm std}$.}

If this hypothesis  is satisfied, we say that $\gg, \G$ is  {\bf $I$-central}.

\noindent
{\bf THE MAIN THEOREM 1.3.}  {\sl
If $\gg, \G$ is an $I$-central  \Ga-Dirichlet polynomial on $\Sn$, then the determinant is majorized:}
$$
\gg(A)^{1\over N} \ \geq\ \gg(I)^{1\over N} (\det\, A)^{1\over n} \qquad \forall\, A>0.
\eqno{(1.1)}
$$

Of course, the constant $\gg(I)^{1\over N}$ is the best possible.

\noindent
{\bf Note. 1.4.}  In our former paper \HHHL\ we needed that the polynomial we were considering
had all of its coefficients positive.   This is not true of even the simplest \Ga-Dirichlet polynomial,
the determinant of a $2\times 2$-matrix: $\det A = a_{11}a_{22}-a_{12}^2$.  As a result we assumed
orthogonal or unitary invariance (as in \HHHL), and this allowed us to work in eigenvalue space,
where all the coefficients are $\geq 0$.

\noindent
{\bf Note. 1.5.}
The Main Theorem can be viewed in the following way. In [\HHHHL, Prop. 11.4], the authors show the following enhancement
of Lemma 4 in \GP.

  \Prop{1.6} {\sl  Let $\gg$ be a \Ga-Dirichlet operator of degree $N$  with \Ga\ cone $\G\ss \Sn$.  
  For any constant $\g > 0$,  the following are equivalent:

(1) \ \ (Majorization of the Determinant)
$$
\gg(A)^{1\over N}\ \geq\     n\g^{1\over n} (\det A)^{1 \over n} \qquad\forall\, A>0
$$

(2)\ \ (Lower Bound for the Determinant of the Gradient)
$$
 D_B \gg^{1\over N} > 0 \and     \det \,( D_B \gg^{1\over N} ) \ \geq\ \g \qquad \forall\, B \in \G.
$$
}

\Cor {1.7} {\sl
The conclusion (1.1) of the Main Theorem can be restated in the following equivalent form:}
$$
 D_B \gg^{1\over N} > 0 \and   \det( D_B \gg^{1\over N}) \ \geq\ \g \ \equiv\ {\gg(I)^{n\over N} \over n^n}, \quad \forall\, B\in \G.
\eqno{(1.1)'}
$$

\vskip.3in


\centerline
{\bf 2.\  The Proof.} 

There are two steps.  First we give the  proof of (1.1) when $A$ is diagonal by using  the Basic Lemma 2.1 in \HHHL \
(see also Gurvits \Gu).
Then we use this diagonal case to prove the general result.

This Basic Lemma can be stated in the following slightly simpler and equivalent form.  Let $e\equiv (1,1,...,1)\in \rn$.

\noindent
{\bf Basic Lemma 2.1.} {Suppose $p$ is a homogeneous polynomial of degree $N>0$ on $\rn$.

\noindent
If 

(1) \ the coefficients of $p(x)$ and all non-negative, and

(2)\ (The Central Ray Hypothesis )
$$
D_e p \ =\ ke\quad\text{for some $k>0$}
\eqno{(2.1)}
$$

\noindent
are both  satisfied,  then 
$$
p(x)^{1\over N} \ \geq\ p(e)^{1\over N}(x_1\cdots  x_n)^{1\over n}, \quad\forall \, x_1>0,\ x_2>0, \ ...\ , \  x_n>0.
\eqno{(2.2)}
$$
}

\noindent
{\bf Note.}  By Euler's Formula for homogeneous polynomials we have  $Np(e) = \bra{D_e p}{e}$,
so if the Central Ray Hypothesis (2.1) holds, then $N, k, n$ and $p(e)$ are related  by $$Np(e) = kn.$$

Our two step procedure is formalized in Theorem 2.3 below where   $\gg$  is replaced by any real homogeneous 
polynomial $F$ on $\Sn$.  Let 
$$
\begin{aligned}
p(x) \ &\equiv \ F(X), \quad \\
{\rm where} \qquad
X \ \equiv 
\left(
\begin{matrix}
x_1 & 0 & \cdots \\
0 & x_2  & \cdots \\
\ & \ & \\
0 & \cdots & x_n
\end{matrix}
\right)
\quad  &{\rm and} \ \ x = (x_1, ... , x_n)\in \rn,
\end{aligned}
\eqno{(2.3)}
$$
be the restriction of $F$ to the set of diagonal matrices $\D \ss \Sn$.
Next, for each orthogonal transformation $h\in {\rm O}(n)$, set
$$
F_h(A) \ \equiv \ F(hAh^t) \and
p_h(x) \ \equiv\ F_h(X) \ =\ F(hXh^t)
\eqno{(2.4)}
$$
for $X\in \D$.

Part (c) of the next result provides an equivalent statement of the Central Ray Hypothesis in terms 
of the polynomials $p_h$ for $h\in  {\rm O}(n)$.

\Lemma {2.2}

\hskip 1in (a) \ \ $D_IF_h \ =\ h^t(D_IF)h$

\hskip 1in (b) \ \ $D_IF_h \ =\ kI \quad\iff\quad D_IF \ =\ kI$

\hskip 1in (c) \ \ $D_IF_h \ =\ kI \quad\iff\quad D_e p_h \ =\ ke \quad \forall\, h \in {\rm O}(n).$

\pf For Part (a)
$$
\begin{aligned} 
\bra {D_I F_h}{A} \ &=\ {d\over dt} \biggr|_{t=0} F_h(I+tA) \ =\  {d\over dt} \biggr|_{t=0} F(h(I+tA) h^t) \\
&=\ {d\over dt} \biggr|_{t=0} F( I + thAh^t)  \ =\ \bra{D_IF}{hAh^t} \\
&= \ \bra{h^t(D_IF)h}{A},  \quad \forall\, A\in \Sn
\end{aligned}
$$

Part (b) is immediate from Part (a).

Proof of  $\Rightarrow$ in Part (c):
$$
\bra{D_e p_h}{x} \ \equiv\  {d\over dt} \biggr|_{t=0} p_h(e+tx) \ =\  {d\over dt} \biggr|_{t=0} F_h(I+tX) \ = \ \bra{D_I F_h}{X}
$$
for all $X\in\D$.  Hence, $D_IF = kI  \ \Rightarrow \ D_ep_h = ke$ for all $h\in{\rm O}(n)$ since $\bra{kI}{X}
= \bra{ke}{x}$.

Proof of  $\Leftarrow$ in Part (c):
Conversely, suppose $D_e p_h = ke$ for all $h\in {\rm O}(n)$.
Given $A\in \Sn$, pick $h\in {\rm O}(n)$ so that 
$$
h^t A h \ =\ \Lambda \ \equiv\ 
\left(
\begin{matrix}
\l_1 & 0 & \cdots \\
0 & \l_2  & \cdots \\
\ & \ & \\
0 & \cdots & \l_n
\end{matrix}
\right)
\quad {\rm is\ diagonal.}
\eqno{(2.5)}
$$
Then with $\l = (\l_1, ... , \l_n)$, we have
$$
\begin{aligned}
&\bra {D_IF}{A} \ =\ \bra {D_IF}{h\L h^t} \ = \ \bra{h^t (D_IF) h}{\L} \quad \text{which by (a)} \\
&= \ \bra{D_I F_h}{\L} \ =\  {d\over dt} \biggr|_{t=0}    F_h(I+t\L) \\
&=\  {d\over dt} \biggr|_{t=0} p_h(e+t\l) \ =\ \bra{D_e p_h}{\l}. \\
\end{aligned}
$$
By the right hand side of part (c) of Lemma 2.2  this last term equals
$\bra{ke}{\l} \ =\ \bra{kI}{A}$. \qed

Continuing with the notation above, we have the following.

\Theorem {2.3}
{\sl  
Suppose $F$ is a real homogeneous polynomial of degree $N$ on $\Sn$ which is 

(1) \ \ $I$-central, \ (See Definition 1.2.)

\noindent
and satisfies the following {\bf coefficient condition}:

(2)\ \ For each $h\in {\rm O}(n)$ the restriction $p_h$ of $F_h$ to the diagonal has coefficients $\geq 0$.

\noindent
Then 
$$
F(A)\on  \ \geq\ F(I)\on (\det A)^{1\over n}  \quad \forall\, A>0.
\eqno{(2.6)}
$$
 }

\noindent 
{\bf Note: } Examples show that the polynomials $p_h(x)$ depend on the choice of $F(A)$, not just $p(x)$.
For example, let $F(A) = a_{11}a_{22} - c a_{12}^2$ with $0\leq c <1$.

\pf
Given $A$ choose $h\in {\rm O}(n)$ so that (2.5) is satisfied.  Then $F(A) = F(h\L h^t)  = p_h(\l)$.

Since $F$ is $I$-central, Lemma 2.2 (b) and (c) implies that $p_h$ is $e$-central.
Hypothosis (2) says that  $p_h$ has coefficients $\geq 0$.  
Therefore by the Basic Lemma 2.1 we have
$$
p_h(\l)\on \ \geq\ p_h(e)\on (\l_1 \cdots\l_n)^{1\over n}.
$$
Substituting $p_h(\l) = F(A), p_h(e) = F(I)$ and $\l_1 \cdots \l_n = \det A$ yields the conclusion (2.6).\qed

To prove the Main Theorem 1.3 it remains to show that  for any $I$-central \Ga-Dirichlet operator $F=\gg$, the coefficient  hypothesis (2) of Theorem 2.3 is satisfied. First we note the following.

\Lemma{2.4} 
{\sl
If $\gg$ is  \Ga-Dirichlet  and $h\in {\rm GL}_n(\bbr)$, then $\gg_h(A) \equiv \gg(hAh^t)$ is also  \Ga-Dirichlet. 
}

\pf
The $I$-eigenvalues of $A\in \Sn$ for the homogeneous polynomial $\gg_h$ are computed as follows:
$$
\begin{aligned}
\gg_h(tI+A) \ &=\  \gg(h(tI+A)h^{t}) \   = \ \gg(thh^t + hAh^{t})   \\
&=\  \gg(h h^t) \prod_j(t+\l_j^{hh^t, \gg}(hAh^t))  \\
&=\  \gg_h(I) \prod_j(t+\l_j^{hh^t, \gg}(hAh^t)).
\end{aligned}
\eqno{(2.7)}
$$
Thus, these eigenvalues are equal to the $hh^t$-eigenvalues of $hAh^t$ for $\gg$.
Note that $hh^t>0$ so $hh^t \in \G(\gg)$.  This proves that the eigenvalues $\l^{I, \gg_h}(A)$ are real, so that 
if $\gg$ is a \Ga\ polynomial, then $\gg_h$ is also a \Ga\ polynomial.  Using the fact that the \Ga\ cone is the set on which all $N$ eigenvalues are strictly positive, this proves that 
$$
A\in \G(\gg_h) \ \ \iff \ \  hAh^t \in \G(\gg), \quad {\rm i.e.}\quad \G(\gg_h) = h^{-1}\G(\gg)(h^t)^{-1}
$$
Finally, since $\Int\cp \ss \G(\gg)$ and $h^{-1}(\Int \cp)(h^t)^{-1} = \Int \cp$, this proves that $\Int\cp\ss \G(\gg_h)$.
Therefore, $\gg_h$ is  \Ga-Dirichlet.
\qed

\noindent
{\bf Note.}   See Appendix A for a further discussion, including when the \Ga-Dirichlet operator 
$\gg_h$, produced by the change of variables $h\in {\rm GL}(n, \bbr)$, is $I$-central.

It remains to show that the  coefficient  condition (2) holds for $\gg_h$.
However, as a consequence of Lemma 2.4 the polynomial $\gg_h$ is \Ga-Dirichlet, so 
we only need to show that  condition (2) holds  for all \Ga-Dirichlet polynomials, which is our final lemma.

\Lemma{2.5} 
{\sl
For any \Ga-Dirichlet polynomial $\gg$, the restriction $$p(x) = \gg(X)$$ to the diagonal $X\in \D$,
has coefficients $\geq 0$.
}

\pf
We use  the following  basic construction of 
 \Ga \ along with the fact  that
$$
{1\over \a_1! \cdots \a_n!}\, {\partial^{|\a|}  p \over \partial x_1^{\a_1}\cdots  \partial x_n^{\a_n}}
\eqno{(2.8)}
$$
equals the coefficient $a_\a$ for each multi-index $\a$ of length $|\a| \equiv \sum_k \a_k = N$.


We  recall that the $a$-eigenvalues of $p$ at $x\in \rn$ are given by  the equation $p(x+ta) = p(a) \prod_{k=1}^N(t + \l_k^p(x))$.

\Prop{2.6.\ (\Ga\ \Gar)}  {\sl  
Suppose $p, \G$ is a \Ga \ polynomial on $V=\rn$ of  degree $N\geq 2$ with $\bbr_{>0}^n \ss\G$.   
Then for $a\in \G$ the polynomial
$$
q(x) \ \equiv \ \bra{D_x p}{a} \ =\ {d\over dt}  p(x+ta) \biggr|_{t=0} \ =\ \sum_{j=1}^n a_j {\partial p \over \partial x_j}(x)
$$
is also \Ga, of  degree $N-1$. 
Moreover, the $a$-eigenvalues of $p$ and $q$ at $x\in\rn$ satisfy:
$$
\l_1^p(x), ... , \l_N^p(x) \ > \ 0 \quad \Rightarrow \quad \ \l_1^q(x), ... , \l_{N-1}^q(x) \ >\ 0
\eqno{(2.9)}
$$
which is equivalent to 
$$
\G(q) \supset \G(p)\supset \bbr_{>0}^n.
\eqno{(2.10)}
$$
}

\pf  For   $a\in\G$ we know that $p$ is hyperbolic in the direction $a$, that is, for every $x\in \rn$ the polynomial $t\mapsto p(x+ta)$ has $N$ real roots.
For every $x\in \rn$ the roots of 
$$
q(x+ta) \ =\  {d\over dt}  p(x+ta)
$$
are the critical points of $p(x+ta)$.  Hence (by the Mean Value Theorem), the $a$-eigenvalues of $q$ interlace the 
$a$-eigenvalues of $p$ for each $x\in \rn$.  This proves that $q$ is \Ga\ of degree $N-1$.  Moreover $\G(p) \ss \G(q)$,
since if $x$ has all positive $p$ eigenvalues, the interlacing $q$-eigenvalues at $x$ must also be positive.
\qed

We shall use Proposition 2.6  
to prove that each  $\a^{\rm th}$  partial of $p$ for $|\a| = N$  is $\geq 0$.
The partial derivative ${\partial p\over \partial x_k}(x)$ is just the directional derivative
in the axis direction $e_k$.  
We  approximate the $\a$ derivative
$$
\left({\partial\over \partial x_1}\right )^{\a_1} \cdots \left ({\partial\over \partial x_n}\right)^{\a_n} 
=
(\nabla_{e_1})^{\a_1} \cdots (\nabla_{e_n})^{\a_n}
$$
with the directional derivative in \Ga\ directions $e_j+\e e \in \G$
$$
(\nabla_{e_1 +\e e})^{\a_1} \cdots (\nabla_{e_n+\e e})^{\a_n}
$$
for $\e>0$.  We now take these derivatives in succession, that is, one at a time.
Since every $e_k+\e e$ lies in the  open set $\bbr^n_{>0}$,  Proposition 2.6 says that each time we take a derivative,
we get another \Ga\  polynomial, of exactly one degree less, whose \Ga\ cone contains the one we just had.
In particular it contains $\bbr^n_{>0}$. After $N-1$ steps we have a linear \Ga\ polynomial
$h(x) \equiv b\cdot x, \  b\in \rn$ with $b\neq 0$, such that the half-space $\G(h)\supset \rn_{>0}$. 
In particular, the inner product $\bra b {e_j + \e e} >0$ for each $j$.  Letting $\e\to 0$ we get that 
the $j^{\rm th}$ component  of $b$, $b_j = \bra b {e_j}\geq 0$ for all $j$.
 Therefore, at the last step, where we take the final derivative, we have a constant which is 
$\geq  0$. This  shows that the $\a^{\rm th}$ partial of $p$ at $x$ is $\geq 0$, and 
Lemma 2.5 is proven.\qed

\noindent
{\bf Note.}  In light of the fact  (proved in Theorem D.1 $(1)'$ of  \HHHL) that the directional derivative of $\gg$ is strictly positive in all directions $B$ in the closed cone $\overline{\G(\gg)}$, except the edge directions, 
one might think that the coefficients $a_\a>0$ for all $|\a|=N$ unless one of the axes belongs to the edge. However this is not true, as the \Ga \ polynomial $\gg(x) = x_1(\e x_1 + x_2)$ in $\bbr^2$ shows.

\noindent
{\bf Proof of the MAIN THEOREM 1.3.} Lemma 2.5 shows that for each of the  \Ga-Dirichlet operators  $\gg_h$,
the restriction  of  $\gg_h$  to the diagonal $\D$ has coefficients $\geq 0$.  That is, the coefficient hypothesis (2) of Theorem 2.3 is satisfied, 
and so Theorem 1.3 holds.\qed

\vskip.3in

\vskip .3in
\def\D{\Delta}

\centerline
{\bf 3. Examples.} 

In this section we examine the family of operators to which our theorems apply. 
In  the restricted cases, considered  in \HHHL, the main result already covered a wide universe of examples.
This was due, in part, because products and radial derivatives 
 take examples to new examples (See Lemmas 7.1 and  Proposition 3.2  in \HHHL).  There is a parallel  situation   here
 without the invariance.

\Prop{3.1} {\sl

(1)\ \ The product of two $I$-central \GD\ operators    on $\Sn$ is again an $I$-central \GD\ operator.

(2)\ \ The redial derivative (or relaxation) of a \GD\ operator on $\Sn$ is also an    $I$-central \GD\ operator.
}

Note that in (1) the degree goes up, whereas in (2) it goes down.  Hence  one can take long sequences of  these operations 
with many "ups and downs" of the degree.  This leads to a large family of operations to which Theorem1.3 applies.

There is a third, very general way of producing examples which includes both constructions in the two steps of our proof of the
Main Theorem 2.3.  This is a third part of Proposition 3.1 which is presented at the end of this section \S 3.

\def\hh{\mathfrak h}
\noindent
{\bf Proof of Proposition 3.1 (1).}
 Suppose $\gg$ and $\hh$ are $I$-hyperbolic on $\Sn$, that is $\gg(tI+A)$ and $\hh(tI+A)$ have only real roots.  Since the roots
 of the product $\gg\hh$ is the union of  the roots of $\gg$ and $\hh$, we see that $\gg\hh$ is also $I$-hyperbolic.

Recall the \Ga\ cone  $\G^\gg$ of $\gg$ is the connected component of $\{\gg>0\}$ containing $I$, and the Dirichlet condition is that
$\cp\ss\G^\gg$.  Of course, the same holds for $\hh$, so clearly $\cp \ss \G^\gg\cap\G^\hh = \G^{\gg\hh}$ since the \Ga\ cones  are convex.   Hence, the product $\gg\hh$
satisfies the Dirichlet condition and is therefore \GD. 

It remains to prove that $\gg\hh$ is $I$-central.
As discussed in \S  1,  this condition has several important equivalent formulations.  One such  is that the {\bf \Ga\ Laplacian}
$
\D^\gg A \ \equiv \ \sum_{j=1}^N \l^\gg_j (A)
$   
and the standard Laplacian 
$
\D A \equiv \sum_{j=1}^N \l_j (A)
$
agree up to a positive factor $\bar k>0$, i.e., $\D^\gg = \bar  k \D$. 
(See Appendix E in \HHHL \  for a complete discussion of the Central Ray Hypothesis 
including several other  equivalent ways of stating it.)  Since $\gg$ and $\hh$ are $I$-central, we have
$\D^\gg = \bar k \D $ and $\D^\hh = \bar \ell \D$ for positive constants $\bar k, \bar \ell$.  Since the set of eigenvalues of $\gg\hh$
is exactly the union of the eigenvalues of $\gg$ and $\hh$,  it  follows that
$\D^{\gg\hh} = \D^\gg + \D^\hh  = (\bar k+\bar \ell)\D$, and so $\gg\hh$ is $I$-central.  \qed 

\noindent
{\bf Remark 3.2. (Concerning the Central Ray Hypothesis)}  {\sl The Central Ray Hypothesis is a weak form of symmetry for the operator.} As seen in \HHHL, \GD\ operators invariant by O$(n)$, U$(n/2)$ or Sp$(n/4)$ are $I$-central.  

If $\gg$ is $I$-central, then $\gg$ is complete, that is, it depends on all the variables ($\iff$ it is not pulled back from a proper vector 
subspace).
This is one reason to 
anticipate that $I$-central \GD\ operators have interior regularity.  Certainly, without completeness, interior regularity fails 
utterly.  If, for example, $\gg(A)$ is independent of $A_{1j}$,  $j=1, ... , n$,  then any continuous function of $x_1$ is  $\gg$-harmonic.

\noindent
{\bf Proof of Proposition 3.1 (2).}  \Ga's result that the radial derivative $\gg'(tI+A) = {d\over dt} \gg(tI+A)$ is again \GD\
is proved in Lemma 2.6.  
To see that $\gg'$ is $I$-central, we consider the family of radial derivatives 
$\gg^{(\ell)}(tI+A) = {d^\ell\over dt^{\ell}}\gg(tI+A)$. 
Following \Gar\ and \Ren, the roots of $\gg^{(\ell+1)}(tI+A)$ are the critical points of $\gg^{(\ell)}(tI+A)$
 and hence intertwine  the roots of $\gg^{(\ell)}(tI+A)$.  
Hence, the associated \Ga\ cones are nested: 
$\G^\gg \ss\G^{\gg'} \ss  \G^{\gg''} \ss\cdots\ss \G^{\gg^{(N-1)}}= H$,  a half space.  Since $\gg$ is $I$-central,
$H = \{A : \D(A)\geq 0\}$.  Now for any $\ell$ we have that $(\gg^{(\ell)})^{(N-\ell-1)} = \gg^{(N-1)}$, from which 
it follows that $\D^{\gg^{(\ell)}} = \bar k \D$, for some $\bar k>0$.\qed

\Remark{3.3. Comments on the Coefficient Condition (2) of Theorem 2.3} 
There are plenty of examples of polynomial operators $F$ satisfying the hypotheses (1) and (2) of Theorem 2.3, and therefore
the determinant majorization (2.6), but which are NOT \Ga.  
For instance 
$$
F(A) \ \equiv \  \|A\|^2 \det A
$$
is O$(n)$-invariant, so that $F_h=F$ and hence $p_h(x) = p(x) = |x|^2x_1\cdots x_n$ has coefficients $\geq 0$ for all
$h\in {\rm O}(n)$. 
However, it remains to determine how interesting they are.

We now want to discuss the following specific example.

\Ex {3.4. (The Lagrangian Monge-Amp\`ere Operator)}
This is an \Ga-Dirichlet operator MA$_{\rm LAG}$ of degree $N=2^n$ on $\bbr^{2n} = \bbc^n$,
which was introduced in \HHLL.   It makes sense on any symplectic manifold with a Gromov metric.

There are two reasons for updating   this particular case.   A determinantal  inequality for this operator was
done in Proposition 8.1 in \HHHL\ with a proof independent of the main theorem in that paper.   Here we show that it is a direct 
consequence of Theorem 1.3 above.
Moreover,  the inequality established in \HHHL\  for MA$_{\rm LAG}$ was not sharp, while the one
given by Theorem 3.1 is.

In \HHHL\ we proved the determinantal inequality
$$
{\rm MA}_{\rm LAG}(A)^{1\over N} \ \geq\ \det(A)^{1\over 2n} \qquad \forall\, A>0.
\eqno{(3.1)}
$$
It is easy to compute that ${\rm MA}_{\rm LAG}(I)^{1\over N} = n$

\Lemma{3.5}
{\sl
The operator ${\rm MA}_{\rm LAG}$ is $I$-central, and hence}
$$
{\rm MA}_{\rm LAG}(A)^{1\over N} \ \geq\ n\, \det(A)^{1\over 2n} \qquad \forall\, A>0.
\eqno{(3.2)}
$$

This is an improvement of (3.1) and is sharp (one has equality when $A=I$).
By the way the proof of this Lemma is easier than the proof of Prop. 8.1 in \HHHL.

 \pf
 The $I$-eigenvalues of $A\in \Sn$ for ${\rm MA}_{\rm LAG}$  were computed in Lemma 3.7 on \HHLL, and they are
 $$
 \L_{\pm \pm \cdots\pm}(A) \ =\ {\tr\, A \over 2} \pm \l_1\pm \l_2 \pm \cdots \pm \l_n
 \eqno{(3.3)}
$$
where $\l_1, ... , \l_n$ are the non-negative eigenvalues  of the skew-hermitian part of $A$.
If we sum (3.3) over the $2^n$ choices of the $\pm$ signs, the $\l_j$'s drop out, and we have
$$
\sum_{\pm}  \L_{\pm \pm \cdots\pm}(A) \ =\ 2^{n-1} \tr\, A.
 \eqno{(3.4)}
$$
So the \Ga\ Laplacian for ${\rm MA}_{\rm LAG}$ is just $2^{n-1}$ times the standard Laplacian of $A$.
This is one of the equivalent ways of saying that ${\rm MA}_{\rm LAG}$  is $I$-central
(See Definition 1.2.).\qed



Finally we add a  third method of constructing examples, which will give a completion to Proposition 3.2.
The two methods used in the proof of the Main Theorem to construct \Ga-Dirichlet operators
(Lemma 2.4 and Proposition 2.6)
can be combined by  replacing  the maps $\Sn \to \bbd$  and $A\mapsto hAh^t$, by a general linear transformation 
$L:\Sn \to \Sn$.  This gives us a third method of constructing $I$-central \GD \ operators, which we address as
a third conclusion (3)  to Proposition 3.1.

\Prop{3.1}  {\sl Under the assumptions of Proposition 3.1 one has the following.

\noindent
(3)\ \  Suppose $L\in \End(\Sn)$ and $\gg$ is a homogeneous polynomial on $\Sn$ of degree $N$.
Consider the transform $\gg_L(A) \equiv \gg(L(A))$, which is also homogeneous of degree $N$.   Then

(a) \ \ $\gg_L$ is $B$-hyperbolic $\iff$ $\gg$ is $L(B)$-hyperbolic (for $B\in \Sn$), in which case the $B$-eigenvalues are given by

(b) \ \  $\l_j^{B, \gg_L}(A) = \l_j^{L(B), \gg}(L(A)), \ \ j=1, ... ,N$, and the \Ga\ cone is given by

(c) \ \ $\G^{B, \gg_L} = L^{-1}(\G^{L(B), \gg})$, i.e., $A\in \G^{B, \gg_L} \iff L(A) \in \G^{L(B), \gg}$.

The gradient map for $\gg_L$ is given by

(d) \ \ $D_A(\log\, \gg_L) \ =\ L^t(D_{L(A)} \log\, \gg)$.

(e) \ \ If $\gg, \G(\gg)$ is \Ga-Dirichlet and $L(\cp) \ss\overline{\G(\gg)}$, then $\gg_L, \G(\gg_L)$ where 
$\G(\gg_L) =   \G^{I, \gg_L}(\gg_L) = L^{-1}(\G(\gg))$
is also \Ga-Dirichlet.

Finally, If $L(I) = k I$ for some $k>0$, 
and if $\tr A=0 \Rightarrow \tr(L(A))=0$, then   

(f) \ \ $\gg$ is $I$-central \ $\Rightarrow$\  \ $\gg_L$ is $I$-central.   
}



  
  \noindent
  {\bf Proof of (a) and (b).}
  Note that 
  $$
  \gg_L(tB+A) \ =\ \gg_L(B)\prod_{j=1}^N (t+\l_j^{B, \gg_L}(A)), \qquad{\rm where}
  $$
   $$
  \gg_L(tB+A) \ \equiv \  \gg(tL(B) + L(A))\ =\ \gg(L(B))\prod_{j=1}^N (t+\l_j^{L(B), \gg}(L(A)))
  $$ 
 just as in the proof of Lemma 2.4.
  This proves (a) and (b). \qed

  \noindent
  {\bf Proof of (c).}  One has that $\gg_L(B) >0 \iff \gg(L(B))>0$ and 
  $$
  \begin{aligned}
A\in \G^{B, \gg_L} \  &\iff  \ \gg_L(A)>0 \ \ {\rm and}\ \  \l_j^{B, \gg_L}(A)>0,\ \ j=1, ... , N    \\
  &\iff  \ \ \gg(L(A)) >0  \ \ {\rm and}\ \   \l_j^{L(B), \gg}(L(A))>0,\ \ j=1, ... , N    \\
  &\iff \ \ L(A) \in \G^{L(B), \gg}. \qquad\mathqed
  \end{aligned}
  $$

  \noindent
  {\bf Proof of (d).}  
  Computing the chain rule for $\Sn  \overset{L}  \arr \Sn \overset {\log\, \gg} \arr \bbr$
  at $A\in \Sn$ one has
  $$
  \bra{D_A\log\, \gg_L}{B} \ =\ \bra{D_{L(A)}\log\, \gg}{L(B)} \ =\ \bra{L^tD_{L(A)} \log\, \gg}{B},
  $$
  for each $B\in\Sn$.\qed

  \noindent
  {\bf Proof of (e).}  
  Given $P>0$, we must show that $P\in \G(\gg_L) = L^{-1}(\G(\gg))$, i.e., $L(P)\in \G(\gg)$. 
  We are assuming $L(P)>0$, and since $\gg, \G(\gg)$ is \Ga-Dirichlet, we have $L(P)\in \G(\gg)$.\qed

  \noindent
  {\bf Proof of (f).}  
  By (d) with $A=I$, we have
  $$
  D_I\log\, \gg_L \ =\ L^t(D_{L(I)}\log\,\gg) \ =\ L^t(D_{kI}\log\,\gg) \ =\ {1\over k} L^t(D_{I}\log\, \gg).
  $$
  Since $I$ is the central ray for $\gg$, this equals ${N\over kn} L^t(I)$. 
  
  It remains to show that $L^t(I) = kI$.  Suppose $A\in\Sn$, with trace-free part $A_0$, i.e., $A \equiv {\tr A\over n}I + A_0$.
  By hypothesis, $\tr A_0 = 0 \ \Rightarrow\ \tr L(A_0)= \bra I {L(A_0)}=0$.
  Therefore, $\bra {L^t(I)}{A} = \bra{I}{L(A)} = \bra{I}{k{\tr A \over n} I} \ =\ \bra{kI}{A}$ for all $A\in\Sn$.\qed

\vskip .3in


\centerline{\bf Bibliography}

 \noindent
\AO\ \  S. Abja and G. Olive, {\sl Local regularity for concave homogeneous complex
degenerate elliptic equations dominating  the
Monge-Amp\`ere equation},   Ann.\ Mat.\ Pura Appl.\ (4) {\bf201} (2022) no.\ 2, 561-587.

\noindent
\ADO \  S. Abja, S. Dinew and G. Olive, {\sl Uniform estimates for concave
homogeneous complex degenerate elliptic equations comparable to the Monge-Amp\`ere
equation}, Potential Anal. (2022). \newline https://link.springer.com/article/10.1007/s11118-022-10009-w

\noindent
\Gar\    L. G\aa rding, {\sl  Linear hyperbolic differential equations with constant coefficients},
Acta. Math. {\bf 85}   (1951),  2-62.

\noindent
\GP\  B.Guo and D. H.Phong, {\sl On  $L^\infty$-estimates for fully nonlinear partial differential equations}, 
  arXiv:2204.12549v1.

\noindent
\GGP  \ \------------ , {\sl Uniform entropy and energy bounds for fully non-linear equations}, ArXiv:2207.08983v1, July 2022.

\noindent
\GGGP \ \------------ , {\sl Auxiliary Monge-Ampere equations in geometric analysis}, arXiv:2210.13308.

\noindent
\GGGGP \ \------------ , {\sl Uniform $L^{\infty}$ estimates: subequations to fully nonlinear partial differential equations}, arXiv:2401.11572v1

\noindent
\GPT\ B. Guo, D.H. Phong, and F. Tong, {\sl On $L^\infty$-estimates for complex Monge-Amp`ere equations},
Ann. of Math. (2)  {\bf 198} (2023), no.1, 393 - 418.
arXiv:2106.02224.

\noindent
\Gu   \    L. Gurvits, {\sl Van der Waerden/Schrijver-Valiant like conjectures
and stable (aka hyperbolic) homogeneous
polynomials :
one theorem for all},  Electron. J. Combin. {\bf 15} (2008), no. 1, Research Paper 66, 26 pp.

\noindent
{\HL}\ F. Reese Harvey and H. Blaine Lawson, Jr.
 {\sl   Dirichlet duality and the non-linear Dirichlet problem on Riemannian manifolds}, 
  J. Diff. Geom.  {\bf 88} No. 3 (2011), 395-482.  ArXiv:0907.1981.

\noindent
 {\HHL} \  -----------,     {\sl  G\aa rding's theory of hyperbolic polynomials},
   {\sl Communications in Pure and Applied Mathematics}  {\bf 66} no. 7 (2013), 1102-1128.

\noindent
 {\HHLL} \  -----------,      {\sl  Lagrangian potential theory and  a  Lagrangian  equation  of  Monge-Amp`ere  type},  
  pp.\  217-258 in Surveys in Differential Geometry,  edited by H.-D. Cao, J. Li, R. Schoen and S.-T. Yau, vol. 22,
International Press, Sommerville, MA, 2018.     ArXiv:1712.03525

\noindent
{\HHHL}  \  -----------,     {\sl  Determinant majorization and the work of Guo-Phong-Tong and Abja-Olive},
   Calc. of Var. and PDE's  {\bf 62} article number: 153 (2023). ArXiv:2207.01729

\noindent
{{\HHHHL}}  \  -----------,    {\sl Duality for \Ga\  operators and a new inequality of GM/AM type}, (to appear).

\noindent
\KLW \ N. Klemyatin, S. Liang and C. Wang, {\sl On uniform estimates for $n-1)$-form fully nonlinear partial differential equations on compact Hermitian manifolds}, arXiv:2211.13798v2.

\noindent
\Ren \   J.  Renegar, {\sl Hyperbolic programs, and their derivative relaxations},
Found.  Comput.  Math. {\bf 6} (2006) no. 1, 59-79.


\end{document}